\documentclass[10pt]{amsart}

\usepackage{amsfonts,amsmath,amsthm}

\newcommand{\Dcal}{{\mathcal D}}
\newcommand{\Rcal}{{\mathcal R}}

\DeclareMathOperator{\adj}{adj}
\DeclareMathOperator{\sgn}{sgn}
\DeclareMathOperator{\im}{im}
\DeclareMathOperator{\Mat}{M}

\newtheorem*{thm}{Theorem}
\newtheorem{cor}{Corollary}

\begin{document}

\begin{abstract}
Given rings $R\subseteq S$, consider the division closure $\Dcal(R,S)$ and the rational closure
$\Rcal(R,S)$ of $R$ in $S$. If $S$ is commutative, then $\Dcal(R,S)=\Rcal(R,S)=RT^{-1}$, where
$T=\{t\in R \; | \; t^{-1}\in S\}$. We show that this is also true if we assume only that $R$ is
commutative.
\end{abstract}

\author{Mark Grinshpon}
\title{Invertibility of Matrices over Subrings}
\address{Department of Mathematics, Virginia Tech}
\email{mgrinshp@math.vt.edu}
\date{August 13, 2006}
\maketitle

\section*{Introduction}

Let $R\subseteq S$ be rings. The division closure $\Dcal(R,S)$ of $R$ in $S$ is the smallest
subring $D$ of $S$ containing $R$ with the property: if $d\in D$, $d^{-1}\in S$, then $d^{-1}\in
D$. The rational closure $\Rcal(R,S)$ of $R$ in $S$ is the smallest subring $D$ of $S$ containing
$R$ with the property: if $A$ is a matrix over $R$ invertible over $S$, then $A^{-1}$ has all
entries in $D$.

An alternative, and more convenient for our purposes, description of the rational closure is given
by $\Rcal(R,S)=\left\{t\in S\,:\,t\mbox{ appears in }A^{-1}\mbox{ for some matrix }A\mbox{ over
}R\right\}$. It follows from Proposition~7.1.1 and Theorem~7.1.2 in \cite{Cohn} that the set thus
defined is in fact a ring; and by Proposition~3.3 in \cite{Linnell}, this is equivalent to the
definition above.

If $S$ is commutative, then
\begin{equation} \label{allsame} \Dcal(R,S)=\Rcal(R,S)=RT^{-1}, \end{equation}
where $T=\{t\in R \; | \; t^{-1}\in S\}$ and $RT^{-1}=\{rt^{-1} \; | \; r\in R,t\in T \}$, due to
the standard formula for matrix inverses, viz.
\begin{equation} \label{matinv} A^{-1}=(\det A)^{-1}\cdot\adj A=\adj A\cdot(\det A)^{-1}. \end{equation}

Note that in general $RT^{-1}\subseteq\Dcal(R,S)\subseteq\Rcal(R,S)$, and both inclusions can be
proper. But if formula~(\ref{matinv}) holds for all invertible matrices (as is the case for
commutative rings), then $\det A\in T$ and every entry of $A^{-1}$ is in $RT^{-1}$; therefore
$\Rcal(R,S)\subseteq RT^{-1}$, proving (\ref{allsame}).

Is (\ref{allsame}) still true if we assume only that $R$ is commutative? In view of the above
argument, this question can be rephrased as follows: for $A\in\Mat_n(R)$, an $n\times n$ matrix
over $R$, invertible in $\Mat_n(S)$, is $\det A$ invertible in $S$? The answer is yes, i.e. the
following results hold.

\begin{thm}
Let $S$ be a ring, $R\subseteq S$ a subring, and assume that $R$ is commutative. If a matrix
$A\in\Mat_n(R)$ is invertible in $\Mat_n(S)$, then $\det A$ is invertible in $S$.
\end{thm}

\begin{cor}
Let $S$ be a ring, $R\subseteq S$ a subring, and assume that $R$ is commutative. Then
$\Dcal(R,S)=\Rcal(R,S)=RT^{-1}$, a commutative subring of $S$, where $T=\{t\in R \; | \; t^{-1}\in
S\}$ and $RT^{-1}=\{rt^{-1} \; | \; r\in R,t\in T \}$.
\end{cor}

More generally, we can consider (see \cite{Cohn}, Chapter 7) a ring homomorphism $f:R\to S$, not
necessarily an imbedding. The division closure $\Dcal^f(R,S)$ is the smallest subring of $S$
containing $\im f$ and closed under taking inverses of elements invertible in $S$, i.e.
$\Dcal^f(R,S)=\Dcal(\im f,S)$; similarly, the rational closure $\Rcal^f(R,S)=\Rcal(\im f,S)$.

\begin{cor}
Let $f:R\to S$ be a ring homomorphism, and assume that $\im f$ is commutative. Then
$\Dcal^f(R,S)=\Rcal^f(R,S)=(\im f)T^{-1}$, a commutative subring of $S$, where $T=\{t\in\im f \; |
\; t^{-1}\in S\}$.
\end{cor}

This note presents a proof of the above stated results.

\section*{Preliminary considerations}

$A$ is invertible means there exists some $B\in\Mat_n(S)$ such that $AB=BA=I$. In the
commutative case, this would imply $\det(A)\det(B)=1$. But in our setting entries of $B$ lie in an
a priori non-commutative ring $S$, so there is no well-defined determinant of $B$. However, by
mimicking a straightforward proof of the Cauchy-Binet formula (see e.g.\ \cite{CB}), of which this
property of determinants is a special case, it is possible to prove that $\det(A)$ is invertible,
with the inverse given by a ``$\det(B)$'' --- a specific expansion of the $n\times n$ determinant in
which all products are taken in an arbitrary but fixed order.

Notation used in the proof: $\sigma=\{i_1,\ldots,i_n\}\in S_n$ means the permutation in $S_n$
acting via $\sigma(t)=i_t$ for $1\le t\le n$.

\section*{Example: case $2\times 2$}

We have that:
$$
AB=
\begin{pmatrix}
a_{11} & a_{12} \\
a_{21} & a_{22}
\end{pmatrix}
\begin{pmatrix}
b_{11} & b_{12} \\
b_{21} & b_{22}
\end{pmatrix}
=
\begin{pmatrix}
a_{11}b_{11}+a_{12}b_{21} & a_{11}b_{12}+a_{12}b_{22} \\
a_{21}b_{11}+a_{22}b_{21} & a_{21}b_{12}+a_{22}b_{22}
\end{pmatrix}
=
\begin{pmatrix}
1 & 0 \\
0 & 1
\end{pmatrix}
.
$$

For convenience, we will use $d_{ij}$ to refer to the entries of the identity matrix. Let us
compute the determinant of this identity matrix written as the product of $A$ and $B$. Of course,
the result will be 1. But since the entries of $B$ are possibly non-commuting, we need to adopt a
certain way of multiplying and expanding expressions involving $b_{ij}$. Each product occurring in
the expansion of the determinant will be multiplied from left to right and in some sense ``from
inside out''.

Let us start off with the product along the main diagonal $d_{11}d_{22}$. Take
$$
1=d_{11}=a_{11}b_{11}+a_{12}b_{21}.
$$
Multiply it by: $a_{21}$ on the left and $b_{12}$ on the right, $a_{22}$ on the left and $b_{22}$
on the right. Get:
\begin{gather*}
a_{21}b_{12}=a_{21}1b_{12}=a_{21}(a_{11}b_{11}+a_{12}b_{21})b_{12}=a_{21}a_{11}b_{11}b_{12}+a_{21}a_{12}b_{21}b_{12}; \\
a_{22}b_{22}=a_{22}1b_{22}=a_{22}(a_{11}b_{11}+a_{12}b_{21})b_{22}=a_{22}a_{11}b_{11}b_{22}+a_{22}a_{12}b_{21}b_{22}.
\end{gather*}
Summing these up, we obtain $1=d_{11}d_{22}$ as:
\begin{equation*}
\begin{split}
1&=a_{21}b_{12}+a_{22}b_{22}=a_{21}1b_{12}+a_{22}1b_{22} \\
 &=a_{21}a_{11}b_{11}b_{12}+a_{21}a_{12}b_{21}b_{12}+a_{22}a_{11}b_{11}b_{22}+a_{22}a_{12}b_{21}b_{22} \\
 &=\sum_{i,j=1}^2 a_{2j}a_{1i}b_{i1}b_{j2}=\sum_{i,j=1}^2 a_{1i}a_{2j}b_{i1}b_{j2},
\end{split}
\end{equation*}
since the entries of $A$ commute with each other.

Next, let us evaluate the product along the other diagonal in a similar fashion. Take
$$
0=d_{21}=a_{21}b_{11}+a_{22}b_{21}.
$$
Multiply it by: $a_{11}$ on the left and $b_{12}$ on the right, $a_{12}$ on the left and $b_{22}$
on the right. Get:
\begin{gather*}
0=a_{11}0b_{12}=a_{11}(a_{21}b_{11}+a_{22}b_{21})b_{12}=a_{11}a_{21}b_{11}b_{12}+a_{11}a_{22}b_{21}b_{12}; \\
0=a_{12}0b_{22}=a_{12}(a_{21}b_{11}+a_{22}b_{21})b_{22}=a_{12}a_{21}b_{11}b_{22}+a_{11}a_{22}b_{21}b_{22}.
\end{gather*}
Summing these up, we obtain $0=d_{21}d_{12}$ as:
\begin{equation*}
\begin{split}
0&=a_{11}0b_{12}+a_{12}0b_{22} \\
 &=a_{11}a_{21}b_{11}b_{12}+a_{11}a_{22}b_{21}b_{12}+a_{12}a_{21}b_{11}b_{22}+a_{11}a_{22}b_{21}b_{22} \\
 &=\sum_{i,j=1}^2 a_{1j}a_{2i}b_{i1}b_{j2}.
\end{split}
\end{equation*}

Now, the determinant of the identity matrix is $1=d_{11}d_{22}-d_{21}d_{12}$ written as:
\begin{equation*}
\begin{split}
1&=\sum_{i,j=1}^2 a_{1i}a_{2j}b_{i1}b_{j2}-\sum_{i,j=1}^2 a_{1j}a_{2i}b_{i1}b_{j2} \\
 &=\sum_{i,j=1}^2 (a_{1i}a_{2j}-a_{1j}a_{2i})b_{i1}b_{j2}=\sum_{i,j=1}^2 \begin{vmatrix}a_{1i}&a_{1j}\\a_{2i}&a_{2j}\end{vmatrix}b_{i1}b_{j2}.
\end{split}
\end{equation*}
Note that $\begin{vmatrix}a_{1i}&a_{1j}\\a_{2i}&a_{2j}\end{vmatrix}$ equals zero if $i=j$, so the
corresponding terms vanish. And when $i$ and $j$ are distinct, this is $\det(A)$ up to the sign.
Thus:

\begin{equation*}
\begin{split}
1&=\sum_{i,j=1}^2 \begin{vmatrix}a_{1i}&a_{1j}\\a_{2i}&a_{2j}\end{vmatrix}b_{i1}b_{j2}=\sum_{\sigma=\{i,j\}\in S_2}\sgn(\sigma)\det(A)b_{i1}b_{j2} \\
 &=\det(A)\cdot\sum_{\sigma=\{i,j\}\in S_2}\sgn(\sigma)b_{i1}b_{j2}.
\end{split}
\end{equation*}

So $\det(A)$ is invertible from the right. Similarly from the other side.

\section*{The general proof}

We have that:
$$
AB=
\begin{pmatrix}
\sum\limits_{k=1}^n a_{ik}b_{kj}
\end{pmatrix}_{i,j=1}^n
=
\begin{pmatrix}
1 & 0 & \ldots & 0 \\
0 & 1 & \ldots & 0 \\
\vdots & \vdots & \ddots & \vdots \\
0 & 0 & \ldots & 1
\end{pmatrix}
.
$$

As in the  sample $2\times 2$ case, let us compute the determinant of this identity matrix written
as the product of $A$ and $B$. We will multiply each term in the expansion of the determinant from left
to right, i.e.
$$
1=\sum_{\sigma=\{i_1,\ldots,i_n\}\in S_n}\sgn(\sigma)d_{i_11}\cdots d_{i_nn}.
$$

Start with
$$
d_{i_11}=\sum_{k_1=1}^n a_{i_1 k_1}b_{k_11},
$$
which is either 0 or 1. In either case,
\begin{multline*}
d_{i_22}=\sum_{k_2=1}^n a_{i_2 k_2}b_{k_22} \Longrightarrow \\
d_{i_11}d_{i_22}=\sum_{k_2=1}^n a_{i_2 k_2}d_{i_11}b_{k_22}=\sum_{k_1,k_2=1}^n a_{i_2 k_2}a_{i_1
k_1}b_{k_11}b_{k_22}.
\end{multline*}
Proceeding in this fashion, we get
$$
d_{i_11}d_{i_22}\cdots d_{i_nn}=\sum_{k_1,k_2,\ldots,k_n=1}^n a_{i_n k_n}\cdots a_{i_2 k_2}a_{i_1
k_1}b_{k_11}b_{k_22}\cdots b_{k_nn}.
$$

Now, the determinant of the identity matrix can be written as:
\begin{equation*}
\begin{split}
1 & =\sum_{\sigma=\{i_1,\ldots,i_n\}\in S_n}\sgn(\sigma)d_{i_11}\cdots d_{i_nn} \\
  & =\sum_{\sigma=\{i_1,\ldots,i_n\}\in S_n}\sgn(\sigma)\sum_{k_1,k_2,\ldots,k_n=1}^n a_{i_n k_n}\cdots a_{i_2 k_2}a_{i_1 k_1}b_{k_11}b_{k_22}\cdots b_{k_nn} \\
  & =\sum_{k_1,k_2,\ldots,k_n=1}^n\left(\sum_{\sigma=\{i_1,\ldots,i_n\}\in S_n}\sgn(\sigma) a_{i_n k_n}\cdots a_{i_2 k_2}a_{i_1 k_1}\right)b_{k_11}b_{k_22}\cdots b_{k_nn} \\
  & =\sum_{k_1,k_2,\ldots,k_n=1}^n\left(\sum_{\sigma=\{i_1,\ldots,i_n\}\in S_n}\sgn(\sigma) a_{i_1 k_1}a_{i_2 k_2}\cdots a_{i_n k_n}\right)b_{k_11}b_{k_22}\cdots b_{k_nn},
\end{split}
\end{equation*}
since the entries of $A$ commute with each other.

Note that the expression in parentheses is precisely the determinant of the matrix whose columns,
say from left to right, are the columns $k_1, k_2,\ldots,k_n$ of the matrix $A$. If not all $k_1,
k_2,\ldots,k_n$ are distinct, such a determinant is zero. And when they are distinct, this is
$\det(A)$ up to the sign. So we can continue:
\begin{equation*}
\begin{split}
1 & =\cdots \\
  & =\sum_{\tau=\{k_1,k_2,\ldots,k_n\}\in S_n}\left(\sum_{\sigma=\{i_1,\ldots,i_n\}\in S_n}\sgn(\sigma) a_{i_1 k_1}a_{i_2 k_2}\cdots a_{i_n k_n}\right)b_{k_11}b_{k_22}\cdots b_{k_nn} \\
  & =\sum_{\tau=\{k_1,k_2,\ldots,k_n\}\in S_n}\sgn(\tau)\det(A)b_{k_11}b_{k_22}\cdots b_{k_nn} \\
  & =\det(A)\cdot\sum_{\tau=\{k_1,k_2,\ldots,k_n\}\in S_n}\sgn(\tau)b_{k_11}b_{k_22}\cdots b_{k_nn}.
\end{split}
\end{equation*}

So $\det(A)$ is invertible from the right. Similarly from the other side.

\section*{Some consequences}

For brevity, set $s=(\det A)^{-1}$. Note that while all entries of $A$ lie in the commutative ring
$R$, and of course so does $\det(A)$, $s$ does not have to be in $R$.

Recall that $A\cdot\adj A=\adj A\cdot A=(\det A)I$, where $\adj A$ is the adjoint matrix of $A$.
Multiplying this from one or the other side by $B$ and then by $s$, we get:
\begin{gather*}
A\cdot \adj A=(\det A)I \quad \Longrightarrow \quad \adj A=B(\det A) \quad \Longrightarrow \quad (\adj A)s=B; \\
\adj A\cdot A=(\det A)I \quad \Longrightarrow \quad \adj A=(\det A)B \quad \Longrightarrow \quad s(\adj A)=B;
\end{gather*}
which is the standard formula for the inverse matrix. This shows that the entries of $B$ lie in
$RT^{-1}$, where $T=\{t\in R \; | \; t^{-1}\in S\}$, and it is easy to see that $RT^{-1}$ is a
commutative subring of $S$.

\section*{Acknowledgements}

I would like to thank Prof.~Peter Linnell for bringing this question to my attention, and for his
help and advice during my work.

\end{document}